\documentclass[11pt,a4paper]{article}
\usepackage{amsmath}
\usepackage{amsfonts}
\usepackage{amsmath,mathtools}
\usepackage{amssymb}
\usepackage{sectsty}
\usepackage{graphicx}
\usepackage{caption}
\usepackage[T1]{fontenc}
\usepackage{bbm}
\usepackage[utf8]{inputenc}
\usepackage{authblk}
\sectionfont{\large}
\numberwithin{equation}{section}
\usepackage[font=scriptsize]{caption}
\usepackage{fancyhdr}
\newcommand*{\TitleFont}{%
      \usefont{\encodingdefault}{\rmdefault}{b}{n}%
      \fontsize{14}{20}%
      \selectfont}

\begin{document}
\title{\TitleFont Shrinkability of Decomposition of $S^n$ Having Arbitrarily Small Neighborhoods with ($n-1$)-Sphere Frontiers }
\author{SHIJIE GU}
\date{}
\maketitle

\begin{abstract}
Let $G$ be a usc decomposition of $S^n$, $H_G$ denote the set of nondegenerate elements and $\pi$ be the natural projection of $S^n$ onto $S^n/G$. Suppose that each point in the decomposition space has arbitrarily small neighborhoods with ($n-1$)-sphere frontiers or boundaries which miss $\pi(H_G)$. If all the arcs are tame in the particular area on the boundary of an $n$-cell $C$ in $S^n$, then this paper shows that this condition implies $S^n/G$ is homeomorphic to $S^n$ ($n\geq 4$). This answers a weak form of a conjecture asked by Daverman [2, p. 61]. In the case of $n=3$, the strong form of the conjecture has an affirmative answer from Woodruff [11].
\end{abstract}

\section{Introduction}
Consider the following open question asked by Daverman [2, p. 61]:

\textbf{Conjecture} (strong): \textit{Suppose} $G$ \textit{is a usc decomposition of} $E^n$ \textit{such that for each} $g\in H_G$ \textit{and each open neighborhood} $W_g$ \textit{of} $g$ \textit{there exists a neighborhood} $U_g$ with $g\subset
U_g\subset W_g$, \textit{where the frontier of} $U_g$ \textit{is an} $(n-1)$-\textit{sphere missing} $N_G$. \textit{G is shrinkable.} 

The difficulty of this conjecture is namely the wildness. In order to avoid the wildness, the conjecture can be modified to a weak form. And that would be our Theorem 1 in section 3.

The technique of proving shrinkability is inspired by Woodruff [11] and Eaton [6]. However, all their methods are only for 3-sphere. Therefore, for $n$-sphere ($n \geq 4$), this technique has to be modified and extended. Our main aim is to generalize the Theorem 2 stated by Woodruff [11] for $n \geq 4$. Therefore, her Lemma 1 and 2 [11] have to be extended. The main part of the proof and extension is section 5. 

In section 6, one will find out that the restriction, all the arcs in Bd $C -(F_1\cup h(F_2))$ are tame, where $C$ is an $n$-cell ($n \geq 4$), $F_1$ and $F_2$ are disjoint 0-dimensional $F_\sigma$-sets in Bd $C$, is necessary if one still wants to apply the useful tool to show shrinkability---disjoint disks (2-cell). Otherwise, there is a counterexample given by Daverman [4], which shows the absence of such kind of disks on the boundary of $C$.

\section{Definitions and Notations}

All decompositions used in this paper are upper semicontinuous ($usc$) defined by Daverman [2]. For a decomposition $G$ of $S^n (n \geq 4)$, the set of nondegenerate elements is denoted by $H_G$, and the natural projection of $S^n$ onto $S^n/G$ by $\pi$. A subset $X\subset S^n$ is $saturated$ (or \textit{G-saturated}) if $\pi^{-1}(\pi(X))=X$. 

The symbol $\rho$ is used to denote the distance fixed on $S^n$ between sets $A$ and $B$ as $\rho(A,B)$; the symbol $\mathbbm{1}$ is the identity map and, for $A \subset X$, $\mathbbm{1}|_A$ is used to denote the inclusion of $A$ in $X$. For $A\subset S^n$ and $\varepsilon>0$, and embedding $h$ of $A$ in $S^n$ is an $\varepsilon$\textit{-homeomorphism} if and only if $\rho(h, \mathbbm{1}|_A)<\varepsilon$. [3]

For $n=3$, each crumpled cube $C$ can be embedded in $S^3$ so that Cl$(S^3-K)$ is a 3-cell have been shown by Hosay [7] and Lininger [9]. In the case of $n\geq 4$, the same claim is proved by Daverman [3].

A $crumpled$ $n$-$cube$ $C$ is a space homeomorphic to the union of an ($n-1$)-sphere in $S^n$ and one of its complementary domains; the subset of $C$ consisting of those points at which $C$ is an $n$-manifold (without boundary) is called the $interior$ $of$ $C$, written as Int $C$, and the subset $C-$Int $C$, which corresponds to the given ($n-1$)-sphere, is called the $boundary$ $of$ $C$, written as Bd $C$. A crumpled $n$-cube $C$ is a $closed$ $n$-$cell$-$complement$ if there exists an embedding $h$ of $C$ in $S^n$ such that $S^n-h$(Int $C)$ is an $n$-cell.

Let $A$ be an annulus bounded by ($n-1$)-spheres $\Sigma_1$ and $\Sigma_2$. A homeomorphism $g$ taking $\Sigma_1$ onto $\Sigma_2$ is called $admissible$ if there exists a homotopy $H:$ $S^{n-1} \times I \rightarrow A$ such that $H(S^{n-1} \times 0) = \Sigma _1$; $H(S^{n-1} \times 1)=\Sigma_2$; and for $x \in S^{n-1}$ if $H(x \times 0)=p \in \Sigma_1$, then $H(x \times 1)=g(p) \in \Sigma_2.$

\section{Shrinkability and Homeomorphism}
Suppose all the arcs in Bd $C -(F_1\cup h(F_2))$ are tame, where $C$ is an $n$-cell, we have below conclusions,

\textbf{Theorem 1.}  \textit{Suppose that for any} $p\in S^n/G$ \textit{and open set U containing p there is an open set V such that} $p\in V \subset U$ \textit{and} Bd $V$ \textit{is a }($n-1$)-\textit{sphere missing} $\pi(H_G)$. \textit{Then} $G$ \textit{is shrinkable.} \textit{Moreover,} $S^n/G$ \textit{is homeomorphic to} $S^n$.

The first part of this theorem is an affirmative answer to the weak form (all the arcs in Bd $C -(F_1\cup h(F_2))$ are tame) of conjecture given by Daverman [4]. 
The theorem can be implied by below theorem,

\textbf{Theorem 2.} \textit{Suppose that for any} $p\in \pi(H_G)$ \textit{and open set U containing p there is an open set V such that} $p\in V \subset U$ \textit{and} Bd $V$ \textit{is a }($n-1$)-\textit{sphere missing} $\pi(H_G)$. \textit{Then} $G$ \textit{is shrinkable.} \textit{Moreover,} $S^n/G$ \textit{is homeomorphic to} $S^n$.

The existence of a homeomorphism between $S^n$ and $S^n/G$ follows from Lemma 1 below, which is modified from Brown's Theorems [1], Daverman's Theorem 5.2, Proposition 6.1, 6.2 and 6.5 [2]:

\textbf{Lemma 1.} \textit{Let G be a usc decomposition of the n-sphere} $S^n$.\textit{The following conditions are equivalent:}

(1) \textit{G is shrinkable};

(2) $\pi: S^n\rightarrow S^n/G$ \textit{is a near homeomorphism};

(3) \textit{If $G$ is a cell-like decomposition, then} $S^n/G$ \textit{is homeomorphic to} $S^n$.

\section{Extension of Woodruff's Lemma and Proof of Theorem 2}
In this section, Theorem 2 will be reduced to Lemma 2, which is similar to Woodruff's Lemma 1 ($n = 3$), but generalized to $n$-sphere.
The proof of the Lemma 2 will be given in section 5.

\textbf{Lemma 2.} \textit{Suppose G is a usc decomposition of} $S^n$, $\varepsilon >0$,\textit{ and C is a crumpled cube in} $S^n$ \textit{with} Bd \textit{C an} ($n-1$)\textit{-sphere which fails to meet any nondegenerate element of G. Then there exists a homeomorphism} $h: S^n \rightarrow S^n$\textit{ such that}

(1) $h|S^n-N(C,\varepsilon)=\mathbbm{1}$;

(2) \textit{if} $g\in G$ \textit{and} $g\subset C$, \textit{then} Diam $h(g)<\varepsilon$, \textit{and}

(3) $if$ $g\in G$, $then$ Diam $h(g)<\varepsilon+$ Diam $g$.
\\*
\\*
Proof of Theorem 2. One should be familiar with the work did by Woodruff [11]. We can imitate the similar strategy used by Woodruff to prove her Theorem 2. That is to find a sequence of saturated open cover $\mathcal{U} \supset \mathcal{U}_1 \supset \mathcal{U}_2 \supset \mathcal{U}_3 \supset \mathcal{U}_4$ of $N_G$, where $N_G$ denotes the union of elements of $H_G$ (called \textit{nondegeneracy}). $\mathcal{U}_4$ is a refinement of $\mathcal{U}_3$ and satifies Woodruff's 5 properties for $\mathcal{U}_4$. Then, in order to apply Lemma 1, one need extend the 2-sphere to ($n-1$)-sphere and $S^3$ to $S^{n}$ in both paragraph 1 and paragraph 5 in Woodruff's proof. 

Define a homeomorphism $\phi$ of $S^n$ to $S^n$. Establish a sequence of homeomorphism $f^i$ fixed off the corresponding $U^{i}_{3}$, where the collection $U^{i}_{3}$'s are separated sets of an open cover $\mathcal{U}_3$. Then, each nondegenerate elements in $U^{i}_{3}$ are covered by $\{U^{1}_{4},\dots,U^{n}_{4}\}\subset \mathcal{U}_4$. The homeomorphism $f^i$ will shrink the diameter of each nondegenerate element to less than a given number. Define $f^i$ as the composition $f^{i}_{n}\dots f^{i}_j\dots f^{i}_{1}$. Each $f^{i}_{j}$ will shrink diameters of elements in $U^{j}_{4}$ to a specified size. The theorem can be proved by repeating such procedure. Define $f^{i}_{j}$. By applying Lemma 2, for each $f^{i}_{j}$ the crumpled cube $X$ in $S^n$ is Cl $f^{i}_{j-1}\dots f^{i}_1(U^{j}_4)$. The following steps are as same as Woodruff's proof.

\begin{flushright}
$\blacksquare$
\end{flushright}

\section{Proof of Lemma 2}
It's not hard to find that the key point to prove Theorem 2 is define $f^{i}_{j}$; in other word, one need to generalize the Lemma 2 for $n \geq 4$. Imitating the proof of the case of $n=3$, we can derive a proof for $n>4$. Although this proof is also strongly inspired by Woodruff's work, there are still some distinctions in details which need be noted.

Proof of Lemma 2. Let $S^n$ be the union of crumpled $n$-cubes $C$ and $C^*=$ Cl $(S^n-C)$. (This can be proved by definition and Jordan-Brouwer Separation Theorem.) Consider the diagram

\begin{figure}[h!]
         \centering
         \includegraphics[width=6cm,height=4cm]{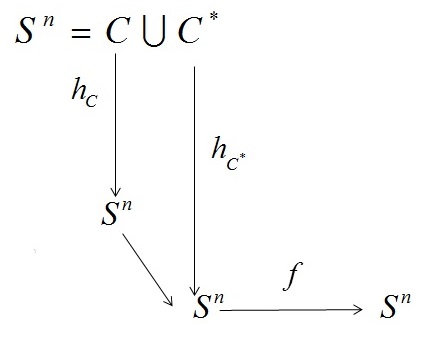}
         \caption{}
         \label{fig:1}
\end{figure}

The map $h_C$ is a reembedding of $C$ in $S^n$ ($n\geq 4$). By Daverman's Reembedding Theorem 6.1 [3], which states: \textit{Let $C$ denote a crumpled n-cube in $S^n$ $(n\geq 4)$. For each $\delta >0$, there exists an embedding $h$ of $C$ in $S^n$ such that $\rho(h,\mathbbm{1}|C) < \delta $ and $S^n-h(\mathrm{Int} C)$ is an n-cell.}

Let $p\in$ Int $C$ and require that $h_C$ is $\mathbbm{1}$ on a small neighborhood of $p$. Let $\theta $ be a homeomorphism of $S^n$ onto itself taking $h_C(C)$ to a set of diameter less than the min\{$\varepsilon/2,$ $\rho(p,$Bd $C)$\} and not moving points in a small neighborhood of $p$. Note that all nondegenerate elements of $\theta h_C(C)$ are Int $C$-small.

Next, we need apply Daverman's Reembedding Theorem to $C^*$ to get the reembedding $h_{C^*}$ of $C^*$ in $S^n$. Choose a $\delta$ which is less than min\{$\varepsilon/2$, $\rho(C^*, \theta h_C(C))$\}. Let $h_{C^*}$ be $\mathbbm{1}$ on $S^n-N(C;\varepsilon)$. The $\varepsilon/2$ condition restrains the growth of nondegenerate elements in $C^*$. Meanwhile, the condition on $\delta$ successfully guarantees $\theta h_C(C) \cap h_{C^*}(C^*)=\emptyset$.

The motion controls imply that $\theta h_C(C)$ and $h_{C^*}(C^*)$ are disjoint crumpled cubes in $S^n$, and the closure of the complement of each is $n$-cell. Denote Cl$(S^n-\theta h_C(C)-h_{C^*}(C^*))$ by $A$. In order to complete the proof of Lemma 2, we have to apply Lemma 3, which is the generalized Woodruff's Lemma 2 [11].
\\*

\textbf{Lemma 3.} \textit{Suppose that in $S^n$ an annulus $A$ is bounded by $(n-1)$-spheres $\Sigma$ and $\varphi(\Sigma)$, where $\varphi$ is an admissible homeomorphism on $\Sigma$; $F_1$ and $F_2$ and disjoint $0$-dimensional $F_\sigma$-sets in $\Sigma$ such that
$F_1 \cup \varphi(F_2) \cup$ $\mathrm{Ext}$ $A$ is $1$-ULC; $U$ is an open set containing $A$; and $\varepsilon >0$. Furthmore, suppose that there exists a decomposition $G(usc)$ of $S^n$ such that}

(1)$f|S^n-U=\mathbbm{1}$,

(2)$f|S^n-A$ \textit{is a homeomorphism onto} $S^n-f(A)$,

(3)$f|\Sigma = f\varphi$ \textit{and} $f\varphi$ \textit{is a homeomorphism onto} $f(A)$ \textit{and}

(4) for $g\in H_G$, Diam $f(g)<\varepsilon+$ Diam $g$.
\\*

Let's show the hypothesis of Lemma 3 can be satisfied. The claim is similar to Woodruff's explanation. But we need modify it. 

In Woodruff's paragraph 1, the 2-sphere $\theta h_X($Bd $X)$ in her Lemma 2 should be extended to ($n-1$)-sphere $\theta h_C($Bd $C)$ 

In paragraph 2, $T$ is changed to be a tame ($n-1$)-sphere in $A$ separating the boundary components. Since Cl $(S^n-\theta h_C(C))$ is $n$-cell. One can consider a homeomorphism $\lambda$ of it onto a polyhedral $n$-cell $P$.
By Generalized Schoenflies Theorem, tame and bicollared are equivalent. So $T$ is bicollared in it. Its image is bicollared in $P$. Then, apply Annulus Theorem, one can show the image is tame. Repeat Annulus Theorem, we
can get two annuli, one is bounded by $\lambda(T)$ and Bd $P$, the other is bounded by $T$ and $\theta h_{C^{*}}($Bd $C^{*})$. It follows that the union of these two annuli whose intersection is an ($n-1$)-sphere belonging to the boundary of each. By Lemma 3.1 [8, p. 167], $A$ is an annulus.

In paragraph 4, Eaton's Mismatch Theorem [6] need be generalized. A sewing $s$ of crumpled $n$-cubes $C$ and $C^*$ has the Mismatch Property if and only if there exist sets $F_{1}'$ and $F_{2}'$ in Bd $C$ such that $F_{1}'$ $\cup$ Int $C$ and $s(F_{2}')$ $\cup$ Int $C^*$ are 1-ULC, and $s(F_0)\cap F_1=\emptyset.$   It follows from Daverman's Theorem 10.1 [10], in case $C$ and $C^*$ of Type 2A [5],  and $n \geq 5$, that is $C \cup_\textrm{Id} C^*=S^n$ if and only if $s$ has the \textit{Mismatch Property}---there exists a 0-dimensional $F_\sigma$ set $F$ in Bd $C$ such that $F$ $\cup$ Int $C$ is 1-ULC, where Id is identity sewing. Hence, $\theta h_X(F_1')$ and $h_{C^*}(F_2')$ are the $F_\sigma$-sets.

These changes complete the generalized hypotheses.
\begin{flushright}
$\blacksquare$
\end{flushright}

Proof of Lemma 3. This is similar to Eaton's proof of his Theorem 4 [6]. That is the repeated application of his Lemma 2. Likewise, apply Woodruff's Lemma 3 and Lemma 4 [11]. However, we need modify Eaton's Lemma 2 or Woodruff's Lemma 3, and add one more hypothesis to produce our new Lemma 4.
\\*

\textbf{Lemma 4}. ($n \geq 5$) \textit{Suppose that $C$ is an $n$-cell in $S^n$ and all the arcs in $\mathrm{Bd}$ $C -(F_1\cup h(F_2))$ are tame. If $D$ is $2$-cell (disk) in Bd $C$, $h$ is a homeomorphism of $2$-cell $D$ onto $\mathrm{Cl}$ $((\mathrm{Bd}$ $C)- D)$ such that $h|\mathrm{Bd}$ $D=\mathbbm{1}$, $F_1$ and $F_2$ are disjoint $0$-dimensional $F_\sigma$-sets in $\mathrm{Int}$ $D$ such that $F_1\cup h(F_2) \cup \mathrm{Ext}$ $C$ is $1$-ULC, $U$ is an open set containing $C-$ $\mathrm{Bd}$ $C$, $\beta >0$, and $\gamma >0$. Furthermore, suppose that there is a decomposition $G$ of $S^n$ such that the nondegenerate elements miss $C$. Then there exists a cellular subdivision $\{D_1, \dots, D_n\}$ of $D$ with mesh less than $\beta$ and a map $f$ of $S^n$ onto $S^n$ such that}

(1) $f|S^n-U=\mathbbm{1}$,

(2) $f|S^n-C$ \textit{is a homeomorphism onto} $S^n-f(C)$,

(3) $both$ $f|D$ $and$ $f|h(D)$ \textit{are homeomorphisms},

(4) $(\bigcup \mathrm{Bd}$ $D_i)\cap (F_1\cup F_2) =\emptyset$,

(5) $f(D)\cup fh(D)=f(\bigcup \mathrm{Bd}$ $D_i)$,

(6) $f|\cup \mathrm{Bd}$ $D_i = fh| \cup \mathrm{Bd}$ $D_i$,

(7) $f(D_i)\cup fh(D_i)$ \textit{bounds a n-cell in} $f(C)$ \textit{of diameter less than} $\beta$, $and$

(8) $for$ $g\in H_G$, Diam $f(g)< \gamma +$ Diam $g$.
\\*

Sketch Proof of Lemma 4. Modify Eaton's proof of his Theorem 4 and Woodruff's proof of her Lemma 3. We will keep all their proofs. However, every time a homeomorphism is performed on $S^n$ not $S^3$. The main strategy used to control the growth of elements are (a) use identity homeomorphism to a neighborhood containing only small elements;(b) each homeomorphism will move points only a small distance: (i) squeeze $C$ into a finite collection of cross-sectionally small cells ($(n-1)$-cells), and then shorten the cells from (i). (ii) Obtain a partial splitting of the cells created in (i). (iii) Shorten and split the cross-sectionally small cells, and repeat this in each finite cells until the whole
$n$-cell is small enough. 

Other parts of the proofs are as same as Woodruff's proof of Lemma 3.
\begin{flushright}
$\blacksquare$
\end{flushright}
Remark 1: One should be careful of the possiblity of wildness. That's why the hypothesis, $C$ is an $n$-cell in $S^n$ and all the arcs in Bd $C -(F_1\cup h(F_2))$ are tame, is necessary. Since it can be shown that each such arc satisfies a 1-LCC condition, i.e., its complement is locally simply 1-connected at each point of the arc. This is to guarantee tameness (or the existence of $2$-cell disks). Otherwise, such tame $2$-cell $D$'s do not always exist. Moreover, $(n-1)$-disks do not exist, either. These can be modified slightly from Daverman's examples to give counterexamples for which there is no small mesh, cellular subdivision $\{D_1,\cdots, D_k\}$  of a certain $(n-1)$-disk  $D$  in  Bd $C$. This shows the straightforward generalization of Eaton's work fails unfortunately.
\\*

In Lemma 3, we suppose that $A$ is an annulus in $S^n$. So we need split or divide $A$ into 2 $n$-cells. In fact, there are two ways to realize it. One is generalize Woodruff's Lemma 4, and that will be our Lemma 5.
The strategy will still work because of the given hypothesis which assures the tameness.

\textbf{Lemma 5}. ($n\geq 4$) \textit{Suppose that} $A$ \textit{is an annulus in} $S^n$, $\Sigma$ \textit{is a} ($n-1$)-\textit{sphere in} Bd $A$, $\varphi$ \textit{is an admissible homeomorphism of} $\Sigma$ \textit{onto} Bd $A-\Sigma$, $F_1$ \textit{and} $F_2$ \textit{are disjoint} 0-\textit{dimensional} $F_{\sigma}$-\textit{sets in} $\Sigma$ such that $F_1 \cup \varphi(F_2) \cup \mathrm{Ext}A$ $is$ $1$-$ULC$, $U$ \textit{is an open set containing $A$}, 
\textit{and} $\gamma>0$. \textit{Furthermore, suppose that there is a decomposition} $G$ \textit{of} $S^n$ \textit{such that the nondegenerate elements miss A}. \textit{Then, there exists an $(n-2)$-sphere} $J$ \textit{standardly embedded in} $\Sigma$ \textit{and a map} $f$ \textit{of} $S^n$ \textit{onto} $S^n$ such that

(1) $f|S^n-U=\mathbbm{1}$,

(2) $f|S^n-A$ \textit{is a homeomorphism onto} $S^n-f(A)$,

(3) \textit{both} $f|\Sigma$ \textit{and} $f|\varphi(\Sigma)$ \textit{are homeomorphisms},

(4) $J \cap (F_1 \cup F_2)=\emptyset$,

(5) $f(\Sigma) \cap f\varphi(\Sigma)=J$,

(6) $f|J=f\varphi|J$,

(7) $f(\Sigma)\cup f\varphi(\Sigma)$ \textit{bounds} 2 $n$-\textit{cells in} $f(A)$, \textit{and}

(8) \textit{for} $g\in H_g$, Diam $f(g)<\gamma +$ Diam $g$. 
\\*

Sketch of the proof. The existence of such $(n-1)$-sphere $J$'s is based on the hypothesis all the arcs in Bd $C -(F_1\cup h(F_2))$ are tame. The strategy of the proof retraces Woodruff's argument (Eaton's Step 1 in his Lemma 2). Step 1, we can have two geometric $n$-cubes which bound an annulus $A'$. By adimissibility, there is a
homeomorphism of such $A'$ into a regular annulus $A$. Then, we can just talk the annulus $A'$ bounded by geometric cubes. Or we can define such a combination of geometric cubes as an \textit{annulus-like n-cube}. Suppose the side of the outer cube is 2, the inner one is 1. Let a $(n-2)$-sphere $J$ separate the boundary of the inner cube into 2 congruent parts.  Consider a distance $t\in [1,2]$, and let $L(t)$ denote the boundary of a geometric cube with side $t$ units. Take a partition of $[1,2]$, $1=t_0<t_1\cdots<t_n=2$. We can have bunch of closed cross-sectional small annulus-like cubes ($(n-1)$-cells) bounded by $L(t_i)$ and $L(t_j)$. Step 2, for each cross-sectional small cell, consider a small neighborhood $V_i$ of the boundary (partial curve Bd $J_i$) which separates annulus-like cube into 2 congruent parts. Push the neighborhood such that it lies in a thin tubular neighborhood. Do it step by step for each partition patiently. We can squeeze the annulus-like cube into 2 $n$-cells. The squeeze is the composition of $\alpha_n\cdots\alpha_1$, where the definition of $\alpha_i$ in details can be referred in Woodruff's paper. Figure 2 represents a schematic diagram of the squeezing process.
\begin{figure}[h!]
         \centering
         \includegraphics[width=10cm,height=8cm]{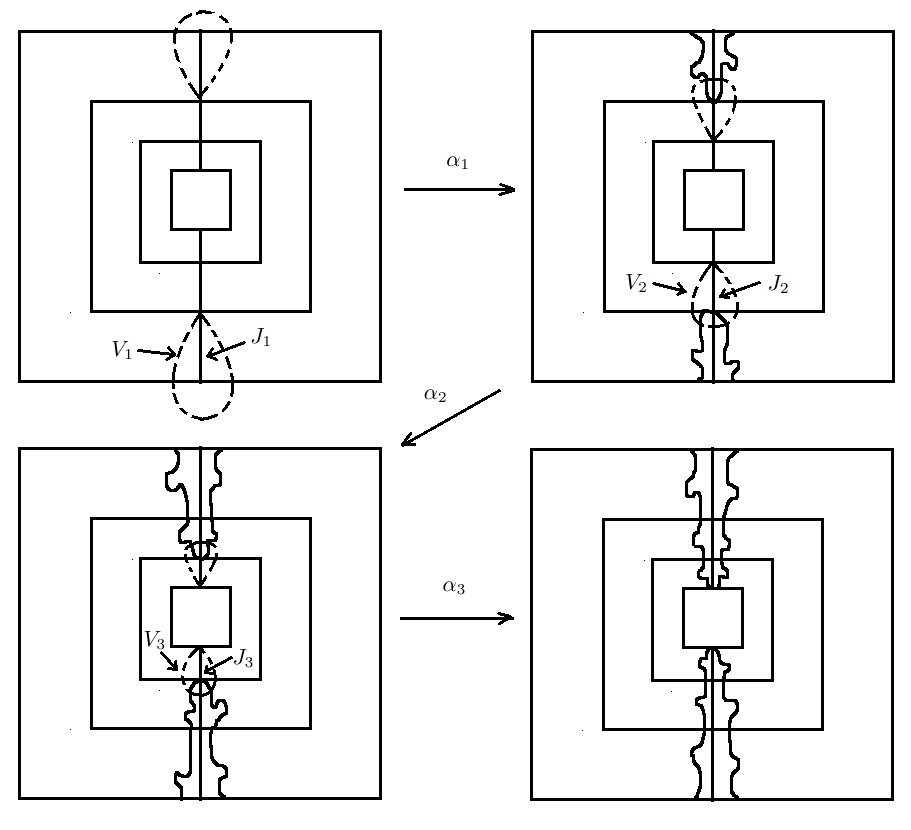}
         \caption{}
         \label{fig:2}
\end{figure}

However, the Lemma 5 can be skipped by using the other approach, which can avoid the wildness, is to consider $n$-dimensional annulus as $S^{n-1} \times [0,1]$. Cover $S^{n-1}$ by two $(n-1)$-cells, $B$ and $B'$, whose intersection is the boundary
of each, i,e.,  $B$ and $B'$ are the upper and lower hemispheres of  $S^{n-1}$. Then $S^{n-1} \times  [0,1]$  
is the union of the $n$-cells  $B \times [0,1]$ and $B' \times [0,1]$. Apply Lemma 4 for each $n$-cell, we can get the proof of Lemma 3 at once.
\begin{flushright}

\end{flushright}
Remark 2: The strategy of Woodruff needs the disks, because the cross-sectional cells ($(n-1)$-cells) are generated from Eaton's disk $D$. 
\section{Absence of Tame Disks in Certain Wild Cells}

We will show the condition $C$ is an $n$-cell in $S^n$ and all the arcs in Bd $C -(F_1\cup h(F_2))$ are tame is necessary in our proofs. Otherwise, Daverman has given following examples [4]:
\begin{itemize}
\item[(1)] For $n \geq 4$, there exists an embedding $f$ of $B^n$ in $E^n$ such that, for each 2-cell $D$ in Bd $B^n$, $E^n-f(D)$ is not simply connected.

\item[(2)] There exists an n-cell $B$ in $E^n$ ($n \geq 4$) such that each 2-cell contained in Bd $B$ is wildly embedded in $E^n$.

\item[(3)] For $2 \leq q <n$ and $n \geq 4$ there exists a $q$-cell $Q$ in $E^n$ such that, for each 2-cell $D$ in $Q$, $E^n-D$ is not simply connected. Hence, each 2-cell in $Q$ is wildly embedded in $E^n$.
\end{itemize}
On the contrary, if the disks have simply connected complement in $E^n$, Seebeck proved that each 2-cell in $E^n$ contains tame arcs [10]. This assures the existence of disks.
\\*
\\*
\textbf{Question.} Can the conclusion still be true if the hypothesis is excluded?

\small{
}

\begin{flushleft}
Department of Mathematics \& Statistics, 
\\*
University of Nevada, Reno
\\*
1664 N. Virginia Street
Reno,  NV  89557-0084
\\*
Email: sgu@unr.edu
\end{flushleft}

\end{document}